\input amstex
\input amsppt.sty   
\magnification=\magstep1
\vsize=20cm  
\hsize=13.5cm  

\def\om{\omega}
\def\row#1#2#3{#1_{#2},\ldots,#1_{#3}}
\def\wedgerow#1#2#3{#1_{#2} \wedge \ldots \wedge #1_{#3}}

\def\Om#1{\Omega^{#1}(M;TM)}
\def\Omh#1{\Omega^{#1}_{\text{hor}}(M;TM)}

\def\Ome#1{\Omega^{#1}(M)}
\def\Omeh#1{\Omega^{#1}_{\text{hor}}(M)}

\def\o{\circ}
\def\X{\frak X}

\def\h{_{\text{hor}}}
\topmatter
\title Graded Derivations \\ of the Algebra of Differential Forms \\
Associated with a Connection \endtitle
\author Peter W. Michor \endauthor
\rightheadtext{Graded derivations associated with connections}
\affil Institut f\"ur Mathematik, Universit\"at Wien,\\
Strudlhofgasse 4, A-1090 Wien, Austria \endaffil
\endtopmatter
\document

\heading Introduction \endheading

The central part of calculus on manifolds is usually the
calculus of differential forms and the best known operators are
exterior derivative, Lie derivatives, pullback and insertion
operators. Differential forms  are a graded commutative algebra
and one may ask for the space of graded derivations of it. It was
described by Fr\"olicher and Nijenhuis in \cite{1}, who found that any
such derivation is the sum of a Lie derivation $\Theta (K)$ and
an insertion operator $i(L)$ for tangent bundle valued
differential forms $K, L \in \Om k$. The Lie derivations give
rise to the famous Fr\"olicher-Nijenhuis bracket, an extension
of the Lie bracket for vector fields to a graded Lie algebra
structure on the space $\Omega (M;TM)$ of vector valued
differential forms. The space of graded derivations is a graded
Lie algebra with the graded commutator as bracket, and this is
the natural living ground for even the usual formulas of
calculus of differential forms. In \cite{8} derivations of even
degree were integrated to local flows of automorphisms of the
algebra of differential forms.

In \cite{6} we have investigated the space of all graded derivations
of the graded $\Omega (M)$-module $\Omega (M;E)$ of all vector
bundle valued differential forms. We found that any such
derivation, if a covariant derivative $\nabla$ is fixed, may
uniquely be written as $\Theta_\nabla (K) + i(L) + \mu(\Xi)$ and
that this space of derivations is a very convenient setup for
covariant derivatives, curvature etc. and that one can get the
characteristic classes of the vector bundle in a very
straightforward and simple manner. But the question arose of how
all these nice formulas may be lifted to the linear frame bundle
of the vector bundle. This paper gives an answer.

In \cite{7} we have shown that differential geometry of principal
bundles carries over nicely to principal bundles with structure
group the diffeomorphism group of a fixed manifold $S$, and that
it may be expressed completely in terms of finite dimensional
manifolds, namely as (generalized) connections on fiber bundles
with standard fiber $S$, where the structure group is the whole
diffeomorphism group. But some of the properties of connections
remain true for still more general situations: in the main part
of this paper a connection will be just a fiber projection onto
a (not necessarily integrable) distribution or sub vector bundle
of the tangent bundle. Here curvature is complemented by
cocurvature and the Bianchi identity still holds. In this
situation we determine the graded Lie algebra of all graded
derivations over the horizontal projection of a connection and we
determine their commutation relations. Finally, for a principal
connection on a principal bundle and the induced connection on
an associated bundle we show how one may pass from one to the other.
The final results relate derivations on vector bundle valued
forms and derivations over the horizontal projection of the
algebra of forms on the principal bundle with values in the
standard vector space.

We use \cite{6} and \cite{7} as standard references: even the notation is
the same. Some formulas of \cite{6, section 1} can also be found in
the original paper \cite{1}.

\heading 1. Connections \endheading

\subheading{1.1} Let $(E,p,M,S)$ be a fiber bundle: So $E$, $M$ and~$S$ are
smooth manifolds, $p : E \to M$ is smooth and each $x
\in M$ has an open neighborhood $U$ such that $p^{-1}(x) =: E
\vert U$ is fiber respectingly diffeomorphic to $U \times S$. $S$
is called the {\it standard fiber}. The {\it vertical bundle}, called
$VE$ or $V (E)$, is the kernel of $Tp : TE \to TM$.
It is a sub vector bundle of $TE \to E$. 

  A {\it connection} on $(E,p,M,S)$ is just a fiber linear
projection $\phi : TE \to VE$; so $\phi$ is a 1-form on
$E$ with values in the vertical bundle, $\phi \in
\Omega^1(E;VE)$. The kernel of $\phi$ is called the {\it horizontal subspace}
$ker \phi$. Clearly for each $u \in E$ the mapping $T_up : ker
\phi \to T_{p(u)}M$ is a linear isomorphism, so $(Tp,
\pi_E) : ker \phi \to TM \times_M E $ is a
diffeomorphism, fiber linear in the first coordinate, whose
inverse $\chi : TM \times_ME \to ker \phi \to
TE$ is called the {\it horizontal lifting}. Clearly the
connection $\phi$ can equivalently be described by giving a
horizontal sub vector bundle of $TE \to E$ or by
specifying the horizontal lifting $\chi$ satisfying
$(Tp,\pi_E)\circ \chi = Id_{TM\times_ME}$. 

  The notion of connection described here is thoroughly treated
in \cite{7}. There one can find parallel transport, 
which is not globally defined in general, holonomy groups,
holonomy Lie algebras, a method for recognizing $G$-connections on
associated fiber bundles, classifying spaces for fiber bundles
with fixed standard fiber $S$ and universal connections for
this. Here we want to treat a more general concept of connection.

\subheading{1.2} Let $M$ be a smooth manifold and let $F$ be a sub vector
bundle of its tangent bundle $TM$. Bear in mind the vertical
bundle over the total space of a principal fiber bundle.

\subheading{Definition} A {\it connection} for $F$ is just a smooth
fiber projection $\phi : TM \to F$ which we view as a
1-form $\phi \in \Omega^1(M;TM)$.

So $\phi_x: T_xM\to F_x$ is a linear projection for all
$ x \in M$. $ker \phi =: H$ is a sub vector bundle of constant
rank of $TM$. We call $F$ the {\it vertical bundle} and $H$ the
{\it horizontal bundle}. $h := Id_{TM} - \phi$ is then the
complementary projection, called the {\it horizontal projection}.
A connection $\phi$ as defined here has been called an almost product
structure by Guggenheim and Spencer.

Let $\Omega_{\text{hor}}(M)$ denote the space of all {\it horizontal}
differential forms $\{ \omega : i_X\omega = 0 \,\text{ for } \,X \in
F\}$. Despite its name, this space depends only on $F$, not on
the choice of the connection $\phi$. Likewise, let $\Omega_{\text{ver}}(M)$
be the space of {\it vertical} differential forms $\{ \omega : 
i_X\omega = 0 \,\text{ for } \,X \in H \}$. This space depends
on the connection. 

\subheading{1.3. Curvature} Let 
$$[\quad,\quad] : \Omega^k(M;TM) \times \Omega^l(M;TM) 
\longrightarrow \Omega^{k+l}(M;TM)$$
be the Fr\"olicher-Nijenhuis-bracket as explained in \cite{6} or
in the original \cite{1}. It induces a graded Lie algebra structure on
$\Omega(M;TM) := \bigoplus \Omega^k(M;TM)$ with one-dimensional
center generated by $Id \in \Omega^1(M;TM)$.  For $K,L \in \Omega^1(M;TM)$ 
we have (see \cite{1} or \cite{6, 1.9})
$$\align [ K, L ] (X,Y) &= [ K(X),L(Y) ] - [ K(Y),L(X) ] \\ 
                        &- L([ K(X),Y ] -[ K(Y),X ] ) \\ 
                        &- K([ L(X),Y ] - [ L(Y),X ] ) \\ 
                        &+ (L \circ K + K \circ L)([X,Y ]).\endalign$$
 From this formula it follows immediately that
$$ [ \phi, \phi ] = [ h,h ] = - [ \phi, h ] = 2(R+\bar R),$$
where $R \in \Omega^2_{\text{hor}}(M;TM)$ and $\bar R \in
\Omega^2_{\text{ver}}(M;TM)$ 
are given by $R(X,Y) = \phi [ hX, hY ]$
and $\bar R(X,Y) = h [ \phi X, \phi Y ]$, respectively.
Thus $R$ is {\it the} obstruction against integrability of the
horizontal bundle $H$; $R$ is called the {\it curvature} of the
connection $\phi$. Likewise $\bar R$ is {\it the} obstruction
against integrability of the vertical bundle $F$; we call $\bar R$
the {\it cocurvature} of $\phi$.

$[ \phi, \phi ] = 2(R+\bar R)$ has been called the torsion of
$\phi$ by Fr\"olicher and Nijenhuis.

\proclaim{1.4. Lemma} (Bianchi-Identity) 
$[ R+\bar R,\phi ] = 0$ and $[R,\phi] = i(R)\bar R +  
i(\bar R)R$ 
\endproclaim 
\demo{ Proof:} We have $2[R+\bar R,\phi ] = [[\phi,\phi],\phi] = 0$ by the
graded Jacobi identity. For the second equation we use the
Fr\"olicher-Nijenhuis operators as explained in \cite{6,section1}. 
$2R = \phi.[\phi,\phi] = i([\phi,\phi])\phi$, and by 
\cite{6,1.10.2}  we have 
$i([\phi,\phi])[\phi, \phi] = \lbrack i([\phi,\phi])\phi, \phi
\rbrack - \lbrack \phi, i([\phi,\phi])\phi \rbrack  +
i([\phi,[\phi,\phi]])\phi + i([\phi,[\phi,\phi]])\phi =
2[i([\phi,\phi])\phi, \phi] = 4 [R,\phi]$.
So $ [R,\phi] = \tfrac14 i([\phi,\phi])[\phi,\phi]  = i(R+\bar R)
(R+ \bar R) = i(R)\bar R + i(\bar R)R$, since R has vertical
values and kills vertical values, so $i(R)R = 0$; likewise for
$\bar R$.  \qed\enddemo

\heading 2. Graded derivations for a connection \endheading

\subheading{2.1.} We begin with some algebraic preliminaries.
Let $A = \bigoplus_{k \in {\Bbb Z}} A_k$ be a graded
commutative algebra and let $I : A \to A $ be an
idempotent homomorphism of graded algebras, so we have $I(A_k)
\subset A_k$, $I(a.b) = I(a).I(b)$ and $I \circ I = I$.  
A linear mapping $D : A \to A$ is called a {\it graded
derivation over $I$ of degree $k$}, if $D : A_q \to
A_{q+k}$ and $D(a.b)=D(a).I(b)+(-1)^{k.\mid a \mid}I(a).D(b)$,
where $\mid a \mid$ denotes the degree of $a$.

\proclaim{Lemma} If $D_k$ and $D_l$ are derivations over $I$ of
degree $k$ and $l$ respectively, and if furthermore $D_k$ and
$D_l$ both commute with I, then the graded commutator $\lbrack
D_k, D_l \rbrack := D_k \circ D_l - (-1)^{k.l} D_l \circ D_k$ is
again a derivation over I of degree $k+l$. 

The space $\text{Der}^I(A) = \bigoplus \text{Der}^I_k(A)$ of
derivations over $I$ which commute with $I$ is a graded Lie
subalgebra of $(End(A),\lbrack \quad,\quad \rbrack\,)$.
\endproclaim
The proof is a straightforward computation.

\subheading{2.2.} Let $M$ be a smooth manifold and let $F$ be a
sub vector bundle of $TM$ as considered in section 1. Let $\phi$
be a connection for $F$ and consider its horizontal projection
$h: TM \to H$. We define $h^*: \Omega  (M) \to  \Omega_{\text{hor}}(M)$
by  $$ (h^*\om)(\row X1p) := \om (\row {hX}1p).$$ Then $h^\ast$ is
a surjective graded algebra homomorphism and 
$h^\ast\vert\Omega_{\text{hor}}(M) = Id$.
Thus $h^\ast: \Omega (M)\to \Omega_{\text{hor}}(M) \to  \Omega
(M)$ is an idempotent 
graded algebra homomorphism.

\subheading{2.3} Let now $D$ be a derivation over $h^\ast$ of $\Omega (M)$
of degree $k$. Then $D\vert\Omega^0(M)$ might be nonzero, e.g.
$h^\ast\circ d$. We call D {\it algebraic}, if $D\vert\Omega^0(M) =
0$. Then $D(f.\om) = 0 + f.D(\om)$, so $D$ is of tensorial
character. Since $D$ is a derivation, it is uniquely determined by
$D\vert\Omega^1(M)$. This is given by a vector bundle homomorphism
$T^*M \to  \Lambda^{k+1}T^*M$, which we view as $K \in
\Gamma(\Lambda^{k+1}\otimes TM) = \Omega^{k+1}(M;TM)$. We write
$D=i^h(K)$ to express this dependence.

\proclaim{Lemma} {\bf 1.} We have  
$$\multline (i^h(K)\om) (\row X1{k+p}) = \\
= \tfrac1{(k+1)!\,
(p-1)!} \sum_\sigma \varepsilon(\sigma).\om(K(\row X{\sigma 1}{\sigma
(k+1)}),\row {hX}{\sigma (k+2)}{\sigma(k+p}),\endmultline$$  
and for any $K \in \Omega^{k+1}(M;TM)$ this formula defines a
derivation over $h^\ast$.

{\bf 2.} We have $[ i^h(K),h^\ast] = i^h(K)\circ h^* - h^*\circ
i^h(K) = 0$ if and only if $h \circ K = K \circ \Lambda^{k+1}h:
\Lambda^{k+1}TM \to  TM$. We write $\Omega^{k+1}(M;TM)^h$ for the
space of all $h$-equivariant forms like that. 

{\bf 3.} For $K \in \Om{k+1}^h$ and $L \in \Om{l+1}^h$ we have: 
$[K,L]^{\wedge,h} := i^h(K)L - (-1)^{kl}\, i^h(L)K$ is an element
of $\Om{k+l+1}^h$ and this bracket is a graded Lie algebra
structure on $\Om{*+1}^h$ such that
$i^h([K,L]^{\wedge,h}) = [i^h(K),i^h(L)]$ in $\text{Der}^h\Omega
(M)$.

{\bf 4.} For $K \in \Om{k+1}$ and the usual insertion operator 
\cite{6, 1.2} we have:
$$\align &h^*\circ i(K) = i^h(K \circ \Lambda^{k+1}h) = h^*\circ i^h(K)\\
&i(K) \circ h^*= i^h(h\circ K) = i^h(K) \circ h^* \\
&[i(K),h^\ast] = [i^h(K),h^\ast] = i^h(h \circ K - K \circ \Lambda^{k+1}h)\\
&h^*\circ i(K) \circ h^*= h^*\circ i^h(K) \circ h^*
    = i^h(h \circ K \circ \Lambda^{k+1}h).\endalign$$
\endproclaim
\demo{Proof} {\bf 1.} We first need the following assertion: For
$\om_j \in \Omega^1(M)$ we have
$$ i^h(K)(\wedgerow \om 1p) = \sum_{j=1}^p (-1)^{(j+1).k}
(\om_1\circ h) \wedge \ldots \wedge (\om_j \circ K) \wedge \ldots \wedge
(\om_p \circ h).\tag1$$
This follows by induction on $p$ from the derivation property.
 From (1) we get
$$\multline \quad i^h(K)(\wedgerow \om 1p)(\row X1{k+p}) = \\
=\sum_{j=1}^p (-1)^{(j-1)k}\bigl((\om_1 \circ h) \wedge
    \ldots \wedge (\om_j \circ K) \wedge \ldots \wedge (\om_p \circ h)
    \bigr) (\row X1{k+p}) = \\
=\sum_j (-1)^{(j-1)k}\tfrac1{(k+1)!}\sum_\sigma
    \varepsilon(\sigma) \om_1(hX_{\sigma 1})\ldots \om_j\bigl( K(\row
    X{\sigma j}{\sigma (j+k)})\bigr) \ldots \\
\dots \om_p(hX_{\sigma (k+p)}).
\endmultline\tag2$$ 
Now we consider the following expression: 
$$\multline \quad\sum_{\pi\in{\Cal S}_{k+p}} \varepsilon
    (\pi). (\wedgerow \om 1p)\bigl( K(\row X{\pi 1}{\pi (k+1)}),\row
    {hX}{\pi (k+2)}{\pi (k+p)}\bigr)=\\
= \sum_{\pi} \varepsilon (\pi)\sum_{\rho \in
    {\Cal S}_p} \varepsilon (\rho)
    \om_{\rho 1}\bigl(K(\row X{\pi 1}{\pi(k+1)})\bigr)
    \om_{\rho 2}(hX_{\pi (k+2)}) \ldots 
    \om_{\rho p}(hX_{\pi (k+p)}),
\endmultline\tag3$$ 
where we sum over $\pi \in {\Cal S}_{k+p}$ and $\rho \in {\Cal 
S}_p$. Reshuffling these permutations one may check that 
$\tfrac1{(k+1)!\,(p-1)!}$ times expression (3) equals expression (2). By
linearity assertion 1 follows.

{\bf 2.} For $\om \in \Omega^1(M)$ we have $\om \circ h \circ K
= i^h(K) h^*\om = h^*i^h(K) \om = h^*(\om \circ K) =
\om \circ K \circ \Lambda^{k+1}h$. 

{\bf 3.} We have $[i^h (K),i^h (L)] = i^h ([K,L]^{\wedge,h})$ for a
unique  $[K,L]^{\wedge,h} \in \Om {k+l+1}^h$ by Lemma 2.1.
For $\om \in \Omega^1(M)$ we get then 
$$\multline \om \circ [K,L]^{\wedge,h} 
    = i^h([K,L]^{\wedge,h})\om =  [i^h(K),i^h(L)]\om =\\
= i^h(K)i^h(L)\om - (-1)^{kl}i^h(L)i^h(K)\om 
    = i^h(K)(\om \circ L)-(-1)^{kl}i^h(L)(\om \circ K).
\endmultline$$ 
So by 1 we get 
$$ \multline i^h(K)(\om \circ L)(\row X1{k+l+1}) = \\
=\tfrac1{(k+1)!\,l!}\sum_{\sigma\in{\Cal S}_{k+l+1}}
    \varepsilon (\sigma)
    (\om \circ L)\bigl( K(\row X{\sigma 1}{\sigma (k+l+1)}),\row {hX}{\sigma
    (k+2)}{\sigma (k+l+1)})\bigr)
\endmultline$$ 
and similarly for the second term.

{\bf 4.} All these mappings are derivations over $h^*$ of
$\Omega (M)$. So it suffices to check that they coincide on
$\Omega^1(M)$ which is easy for the first two assertions. The
latter two assertions are formal consequences thereof. \qed\enddemo

\subheading{2.4} For $K \in \Om k $ we have the {\it Lie
derivation} $\Theta (K) = [i(K),d] \in \text{Der}_k\Omega (M)$, where
$d$ denotes exterior derivative. See \cite{6, 1.3}.

\proclaim {Proposition} {\bf 1.} Let $D$ be a derivation over
$h^*$ of $\Omega (M)$ of degree k. Then there are unique
elements $K \in \Om k$ and $L \in \Om {k+1}$ such that 
$$ D = \Theta (K) \circ h^* + i^h(L).$$
$D$ is algebraic if and only if $K=0$.

{\bf 2.} If $D$ is in $\text{Der}_k^h\Omega (M)$ (so $[D,h^*]=0$) then 
$$D = h^* \circ \Theta (K) \circ h^* + i^h(\tilde L)$$ 
for unique $K \in \Omega^k_{\text{hor}}(M;TM)$ and 
$\tilde L \in \Om {k+1}^h$.
$K$ is the same as in 1.
\endproclaim
Define $\Theta^h(K) := h^* \circ \Theta (K) \circ
h^* \in \text{Der}_k^h\Omega (M)$, then in 2 we can write $D =
\Theta^h(K) + i^h(\tilde L)$. 

\demo{Proof} Let $X_j \in {\Cal X}(M)$ be vector fields and
consider the mapping $ev_{(\row X1k)} \o D : C^\infty (M) =
\Omega^0(M) \to  \Omega^k(M) \to  \Omega^0(M) = C^\infty (M)$, given
by $f \mapsto (Df)(\row X1k)$. This map is a derivation of the
algebra $C^\infty (M)$, since we have $D(f.g)(\row X1k) = (Df.g +
f.Dg)(\row X1k) = (Df)(\row X1k).g + f.(Dg)(\row X1k)$. So it is
given by the action of a unique vector field $K(\row X1k)$,
which clearly is an alternating and $C^\infty (M)$-multilinear
expression of the $X_j$; we may thus view $K$ as an element of
$\Om k$. The defining equation for $K$ is $Df = df \o K$ for $f
\in C^\infty (M)$. Now we consider $D - [i(K),d] \o h^*$. This
is a derivation over $h^*$ and vanishes on $\Omega^0(M)$, since
$[i(K),d] h^* f = [i(K),d] f = df\o K = Df$. It is algebraic and
by 2.3 we have $D  - [i(K),d]\o h^* = i^h(L)$ for some unique
$L \in \Om {k+1}$.

Now suppose that $D$ commutes with $h^*$, i.~e. $D\in \text{
Der}_k^h(M)$. Then for all $f\in C^\infty (M)$ we have $df \o
K\o \Lambda^kh = h^*(df\o K) = h^* Df = Dh^*f = Df = df\o K$, so
$K\o \Lambda^kh = K$ or $h^*K = K$ and $K$ is in $\Omega_{\text{
hor}}^k(M;TM)$. Now let us consider $D - h^* [i(K),d]h^* \in \text{
Der}_k^h\Omega (M)$, which is algebraic, since $h^*[i(K),d]h^*f =
h^*(df\o K) = df\o (h^*K) = df\o K = Df$. So $D- h^*[i(K),d]h^*
= i^h(\tilde L)$ for some unique $\tilde L \in \Om {k+1}^h$ by
2.3 again. \qed\enddemo 

\subheading{2.5.} For the connection $\phi$ (respectively. the horizontal
projection $H$) we define the {\it classical covariant
derivative} $D^h := h^* \o d$. Then clearly $D^h: \Omega (M) \to 
\Omega_{\text{hor}}(M)$ and $D^h$ is a derivation over $h^*$, but
$D^h$ does not commute with $h^*$, so it is not an element of
$\text{Der}^h_1\Omega (M)$. Therefore and guided by 2.4 we
define the (new) {\it covariant derivative} as $d^h := h^* \o d
\o h^* = D^h \o h^* = \Theta^h(Id)$. We will consider $d^h$ as
the most important element in $\text{Der}^h_1\Omega (M)$.

\proclaim{Proposition} \roster
\item"{\bf 1.}" $d^h - D^h = i^h(R)$, 
    where R is the curvature.
\item"{\bf 2.}" $[d,h^*] = \Theta (\phi) \o h^* + i^h(R+\bar R)$.
\item"{\bf 3.}" $d \o h^* - d^h = \Theta (\phi) \o h^* + i^h(\bar R)$.
\item"{\bf 4.}" $D^h \o D^h = i^h(R) \o d$.
\item"{\bf 5.}" $[d^h,d^h] = 2 d^h\o d^h = 2 i^h(R) \o d \o h^* =2 h^* \o
    i(R) \o d \o h^*$.
\item"{\bf 6.}" $D^h\vert\Omega_{\text{hor}}(M) = d^h\vert\Omega_{\text{hor}}(M)$.
\endroster
\endproclaim
Assertion 6 shows that not a lot of differential geometry, on
principal fiber bundles e.g., is
changed if we use $d^h$ instead of $D^h$. In this paper we focus
our attention on $d^h$. 

\demo{Proof} 1 is 2 minus 3. Both sides of 2 are derivations
over $h^*$ and it is straightforward to check that they agree on
$\Omega^0(M) = C^\infty (M)$ and $\Omega^1(M)$. So they agree on
the whole of $\Omega (M)$. The same method proves equation 3.
The rest is easy. \qed\enddemo 

\proclaim{2.6. Theorem} Let $K$ be in $\Om k$. Then we have:
\roster
\item"{\bf 1.}" If $K \in \Om k^h$, then $\Theta^h(K) = [i^h(K),d^h]$. 
\item"{\bf 2.}" $h^*\o\Theta (K) = h^*\o\Theta
    (K\o\Lambda^kh) + (-1)^{k-1}i^h(i^h(R)K)$. 
\item"{\bf 3.}" $\Theta^h(K) = \Theta^h(h^*K) + (-1)^{k-1}i^h(i^h(R)(h\o K))$.
\item"{\bf 4.}" If $K\in\Om k^h$, then
    $\Theta^h(K)=\Theta^h(h^*K)= \Theta^h(h\o K)$.
\endroster
Now let $K\in\Omega_{\text{hor}}^k(M;TM)$. Then we have:
\roster
\item"{\bf 5.}" $h^*\o\Theta (K) - \Theta^h(K) = i^h(\phi\o
    [K,\phi ]\o\Lambda h)$. 
\item"{\bf 6.}" $\Theta (K)\o  h^* - \Theta^h(K) =
    i^h((-1)^ki(\bar R)(h\o K) - h\o [K,\phi ])$. 
\item"{\bf 7.}" $[h^*,\Theta(K)] =  i^h(\phi\o
    [K,\phi]\o\Lambda h + h\o [K,\phi ] -(-1)^ki(\bar R)(h\o K))$. 
\item"{\bf 8.}" Let $K_j\in \Omega_{\text{hor}}^{k_j}(M;TM),
    j=1,2$. Then 
    $$\multline [\Theta^h(K_1),\Theta^h(K_2)] = \Theta^h([K_1,K_2]) - \\
    - i^h\bigl ( h\o [K_1,\phi\o [K_2,\phi ]]\o \Lambda h 
        -(-1)^{k_1k_2}h\o [K_2,\phi\o [K_1,\phi ]]\o \Lambda h\bigr).
    \endmultline$$
\endroster
\endproclaim
\subheading{Note} The third term in 8 should play some r\^ole in the
study of deformations of graded Lie algebras.

\demo{Proof: 1} Plug in the definitions.

{\bf 2.} Check that $h^*\Theta (K) f = h^*\Theta (h^*K)f$ for
$f\in C^\infty (M)$. For $\om  \in \Ome 1$ we get $h^*\Theta
(K)\om = h^*  i^h(h^*K)  d \om -
(-1)^{k-1}  h^*  d  i(K) \om$, using 2.3.4; also we have $h^* \Theta
(h^*K) \om = h^*  i^h(h^*K)  d \om -(-1)^{k-1}h^*  d h^* 
i(K) \om$. Collecting terms and using 2.5.1 one obtains the
result. 

{\bf 3 } follows from 2 and 2.3.4.\quad
{\bf 4 } is easy.

{\bf 5. } It is easily checked that the left hand side is an
algebraic derivation over $h^*$. So it suffices to show that
both sides coincide if applied to any $\om \in \Ome 1$. This
will be done by induction on $k$. The case $k=0$ is easy.
Suppose $k\geq 1$ and choose $X\in\X (M)$. Then we have in turn:
$$\align &i(X)  h^* \Theta(K) \om = i^h(hX) \Theta  (K) \om = h^* 
    i(hX) \Theta (K) \om \\
&\qquad = h^* [i(hX),\Theta (K)] \om + (-1)^k h^* \Theta (K)  i(hX) \om \\
&\qquad = h^*\Theta(i(hX)K)\om +(-1)^kh^*i([hX,K])\om 
    +(-1)^k  h^* \Theta (K)(\om (hX))\\
 &i(X)  h^* \Theta (K)  h^*\om = i(X) h^* \Theta (K)(\om\o h) \\
&\qquad= h^* \Theta(i(hX)K)(\om\o h) + (-1)^kh^*i([hX,K])(\om\o h) 
    +(-1)^k  h^* \Theta (K)(\om (hX))\\
&i(X)\bigl ( h^* \Theta (K)- h^* \Theta (K) h^*\bigr ) \om \\
&\qquad=\bigl ( h^* \Theta (i(hX)K)-h^* \Theta(i(hX)K)h^*\bigr)\om 
    +(-1)^k  h^*  i([hX,K])(\om\o \phi)\\
&\qquad = i^h\bigl ( \phi\o [i(hX)K,\phi]\o \Lambda h\bigr)\om 
    + (-1)^k  h^*  i([hX,K])(\om\o\phi)\qquad
    \text{ by induction on $k$.}\\
&i(X)  i^h\bigl (\phi\o [K,\phi]\o\Lambda h\bigr ) \om 
    = i(X)\bigl (\om\o\phi\o [K,\phi ]\o\Lambda h\bigr ) 
    = (\om\o\phi)\bigl ( i(hX)[K,\phi])\o\Lambda h\bigr )\\
&\qquad =(\om\o\phi )\bigl ([i(hX)K,\phi]+(-1)^k[K,i(hX)\phi =0] 
    - (-1)^ki([K,hX]) \phi  + \\
&\qquad\qquad\qquad\qquad+(-1)^{k-1}i([\phi ,hX])K \bigr )\o\Lambda h, 
    \qquad \text{ by \cite{6, 1.10.2}}\\
&\qquad= i^h\bigl (\phi\o[i(hX)K,\phi]\o\Lambda
    h\bigr ) +(-1)^k   h^* i([hX,K])(\om\o\phi ) +\\
&\qquad\qquad\qquad\qquad+(-1)^k \om\o\phi\o (i([\phi,hX])K)\o\Lambda h.
\endalign$$
Now $[\phi ,hX](Y)=[\phi Y,hX]-\phi ([Y,hX])$ by \cite{6,1.9}.
Therefore $(i([\phi ,hX])K)\o\Lambda h = i(h\o
[\phi,hX]) K\o\Lambda h = 0$, and $i(X)\bigl (h^*\o\Theta (K) -
h^*\o\Theta (K)\o h^*\bigr )\om = i(X)  i^h(\phi\o
[K,\phi]\o\Lambda h) \om$ for all $X\in\X (M)$. Since the $i(X)$
jointly separate points, the equation follows.

{\bf 6.} $K$ horizontal implies that $\Theta (K)\o h^*
-h^*\o\Theta (K)\o h^*$ is algebraic. So again it suffices to
check that both sides of equation 6 agree when applied to an
arbitrary 1-form $\om\in\Ome 1$. This will again be proved by
induction on $k$. We start with the case $k=0$: Let $K=X$ and
$Y\in \X (M)$. 
$$\multline \bigl (\Theta (X)  h^* \om\bigr )(Y) - \bigl
    (h^* \Theta(X) \,h^* \om\bigr )(Y) =\\
= \bigl (\Theta (X)(\om\o h)\bigr )(Y - hY) 
    = X.\om (h\phi Y) - (\om\o h)([X,\phi Y]) =
    0-\om (h[X,\phi Y])
\endmultline$$
Note that $[X,\phi ](Y)=[X,\phi Y]-\phi  [X,Y]$, so
$h [X,\phi ](Y) = h  [X,\phi Y]$, and therefore
$\bigl (i^h(0-h  [X,\phi ]) \om\bigr )(Y)=-\om
(h [X,\phi Y])$.

Now we treat the case $k>0$. For $X\in\X (M)$ we have in turn:
$$\multline i(X) \Theta (K)  h^* \om = \bigl
    ((-1)^k\Theta(K)  i(X)+\Theta (i(X)K)+ (-1)^ki([X,K])\bigr
    )(h^* \om )\quad\text{ by \cite{6,1.6}}\\
=(-1)^k \Theta (K)(\om (hX))+\Theta (i(X)K) h^* \om +(-1)^k \om (h[X,K]).
\endmultline$$
We noticed above that $\Theta(K)  f=h^* \Theta (K)f$ 
for $f\in C^\infty (M)$ since $K$ is horizontal.
$$\multline i(X)  h^* \Theta (K)  h^* \om = i^h(hX) \Theta(K)  h^* \om 
    = h^*  i(hX) \Theta (K)  h^* \om,\quad\text{ by 2.3.4.}\\
=(-1)^k  h^* \Theta (K)  (\om (hX))+  h^* \Theta (i(hX)K)  h^* \om +(-1)^k 
    h^*(\om(h[X,K])),\quad\text{ see above.}
    \endmultline$$
$$\multline i(X)\bigl (\Theta (K)  h^* \om - h^* \Theta (K) 
    h^* \om\bigr )= \\
=\Theta (i(X)K)  h^* \om - h^* \Theta (i(hX)K) h^* \om  
    +(-1)^k\bigl (\om (h[X,K])-h^*(\om\o h\o [hX,K])\bigr )=\\
=-i^h(h[i(X)K,\phi])\om +(-1)^{k-1} i^h(h\o i(\bar R)i(X)K)\om +\\
+(-1)^k\om (h\o [X,K]) - (-1)^k\om (h\o (h^*[hX,K])),
\endmultline$$  
by induction, where we also used $i(X)K=i(hX)K$, which holds for
horizontal $K$. 
$$\multline -i(X)  i^h(h\o [K,\phi ]) \om = -i(X)(\om\o h\o
    [K,\phi ])= \om\o h\o (i(X)[K,\phi ])=\\
=\om\o h\o \bigl ([i(X)K,\phi]+(-1)^k[K,i(X)\phi ]
    -(-1)^ki([K,X])\phi +(-1)^{k-1}i([\phi,X])K\bigr )=\\
=-i^h(h\o [i(X)K,\phi ]) \om +(-1)^k \om (h\o [\phi X,K])+0
    +(-1)^k \om (h\o (i([\phi ,X])K),
\endmultline$$
where we used \cite{6, 1.10.2} and $h\o\phi =0$.
$$\align &i(X)i^h(h\o i(\bar R)K) \om = i(X)(\om\o h\o i(\bar R)K) 
    = \om\o h\o (i(X)i(\bar R)K).\\
&i(X)  (\text{left hand side $-$ right hand side}) \om =\\
&\qquad= \om\bigl ((-1)^{k-1}h\o i(\bar R)\o i(X)\o K+
    (-1)^k  h\o [X,K] -(-1)^k  h\o (h^*[hX,K])-\\
&\qquad\qquad\qquad-(-1)^k  h\o [\phi X,K]- (-1)^k h\o i([\phi ,X])K -
    (-1)^k  h\o (i(X)i(\bar R)K)\bigr )\\
&\qquad=(-1)^k(\om\o h)\bigl(-i([\bar
    R,X]^\wedge)K+ [X-\phi X,K] -h^*[hX,K] - i([\phi ,X])K\bigr )=\\
&\qquad=(-1)^k(\om\o h)\bigl (i(0-i(X)\bar R)K+
[hX,K]+ h^*[hX,K]-i([\phi ,X])K\bigr )\\
&-h \bigl (i(i(X)\bar R)K\bigr )(\row Y1k) =-
    \sum_{j=1}^k h  K(Y_1,\ldots,\bar R(X,Y_j),\ldots,Y_k)=\\
&\qquad=-\sum_{j=1}^k h  K(Y_1,\ldots,[\phi
    X,\phi Y_j],\ldots,Y_k), \text{ since $K$ is horizontal.}\\
&-h  [K,hX]((\row Y1k) = -h  [K(\row Y1k),hX] +
    \sum_{j=1}^k h  K(Y_1,\ldots,[Y_j,hX],\ldots,Y_k).\\
&h(h^*[K,hX])(\row Y1k) = h  [K,hX](\row {hY}1k) =\\
&\qquad=h  [K(\row{hY}1k),hX]-\sum h 
    K(hY_1,\ldots,[hY_j,hX],\ldots ,hY_k) =\\
&\qquad=h  [K(\row Y1k),hX]-\sum h  K(Y_1,\ldots,[hY_j,hX],
    \ldots,Y_k),
\endalign$$
since $K$ is horizontal. We also have
$h[\phi,X](Y)=h([\phi Y,X]-\phi [Y,X])
    =h[\phi Y,X]$, and thus we finally get
$$\align &-h(i([\phi ,X]K)(\row Y1k)=-h(i(h[\phi,X])K)(\row Y1k)=\\
&\qquad=-\sum h  K(Y_1,\ldots,h[\phi,X](Y_j),\ldots,Y_k) =
    -\sum h  K(Y_1,\ldots,[\phi Y_j,X],\ldots,Y_k).
\endalign$$
All these sum to zero and equation 6 follows.

{\bf 7. } This is equation 5 minus equation 6.

{\bf 8.} We compute as follows: 
$$\align &[\Theta^h(K_1),\Theta^h(K_2)] =\\
&\qquad=h^* \Theta (K_1)  h^* \Theta (K_2) 
    h^* - (-1)^{k_1k_2}  h^* \Theta (K_2)  h^* \Theta (K_1)  h^* =\\
&\qquad=\bigl (h^* \Theta (K_1) - h^*  i(\phi\o
    [K_1,\phi ])\bigr )\theta (K_2)  h^* -\\
&\qquad\qquad\qquad - (-1)^{k_1k_2}\bigl
    (h^* \Theta (K_2) - h^*  i(\phi\o [K_2,\phi ])\bigr ) \Theta
    (K_1)  h^* ,\quad\text{ by 5},\\
&\qquad=h^*\bigl (\Theta (K_1) \Theta (K_2)
    -(-1)^{k_1k_2} \Theta (K_2) \Theta (K_1)\bigr) h^* - \\
&\qquad\qquad\qquad - h^*\bigl(i(\phi\o [K_1,\phi ])\Theta (K_2) 
    - (-1)^{k_1k_2}i(\phi\o [K_2,\phi ])\Theta (K_1)\bigr )  h^* =\\
&\qquad=h^* \Theta ([K_1,K_2])  h^* +\text{Remainder} =
    \Theta^h([K_1,K_2]) +\text{Remainder.}
\endalign$$
For the following note that $i(\phi\o [K_1,\phi ]) 
h^* = i^h(h\o \phi\o [K_1,\phi])=0$, by 2.3.4.
$$\align&\text{Remainder}=-h^*  i(\phi\o [K_1,\phi ]) \Theta
    (K_2)  h^* + (-1)^{k_1k_2}i(\phi [K_2,\phi ]) \Theta (K_1)  h^*+\\
&\qquad\qquad+(-1)^{k_1k_2}h^* \Theta (K_2)  i(\phi\o [K_1,\phi
    ])  h^* - h^* \Theta(K_1)  i(\phi\o [K_2,\phi ]) 
    h^*,\quad\text{ which are $0$,}\\
&\qquad=-h^* [i(\phi\o [K_1,\phi ]),\Theta (K_2)]  h^*
    +(-1)^{k_1k_2} h^*  [i(\phi\o [K_2,\phi ]),\Theta (K_1)] h^*=\\
&=h^*\bigl (\Theta (i(\phi\o [K_1,\phi ])K_2=0) + (-1)^{k_2}
    i([\phi\o [K_1,\phi],K_2])\bigr ) h^* +\\
&\qquad\qquad+(-1)^{k_1k_2}h^*\bigl (\Theta (i(\phi\o
    [K_2,\phi])K_1=0) +(-1)^{k_1}i([\phi\o [K_2,\phi ],K_1)\bigr )  h^*=\\
&\qquad=-i^h\bigl (h\o [K_1,\phi\o [K_2,\phi ]]\o\Lambda h
    -(-1)^{k_1k_2}   h\o [K_2,\phi\o [K_1,\phi]]\o\Lambda h\bigr),
    \quad\text{by 2.3.4.}
\endalign$$
\enddemo

\proclaim{2.7. Corollary} 
\roster 
\item $h^*\o\Theta (\phi)=i^h(R)$.
\item $\Theta^h(\phi)=0$.
\item $[h^*,\Theta (h)]=-2i^h(R)-i^h(\bar R)$.
\item $h^*\o\Theta (h) = \Theta^h(h) - 2i^h(R)$.
\item $\Theta^h(h) +i^h(\bar R)= \Theta (h)\o h^*$.
\item $[\Theta^h(h),\Theta^h(h)] = 2\Theta^h(R) = 2h^*\o i(R)\o d\o h^*$.
\item $\Theta^h(\bar R) = 0$.
\endroster
\endproclaim
\demo{Proof} 1 follows from 2.6.2. 2 follows from 2.6.3. 3
follows from 2.6.7. 4 follows from 2.6.5. 5 is 4 minus 3. 6
follows from 2.6.8 and some further computation and 7 is
similar.  \qed\enddemo

\proclaim{2.8. Theorem: 1} For $K\in \Om k $ and $L\in\Om
{l+1}^h$ we have 
$$ [i^h(L),\Theta^h(K)] = \Theta^h(i(L)K)+(-1)k  i^h(h\o
[L,K]\o\Lambda h).$$

{\bf 2. } For $L\in \Om {l+1}^h$ and $K_i\in\Om {k_i}$
we have
$$\multline i^h(L)(h^*  [K_1,K_2]) = h^*  [i(L)K_1,K_2] +
    (-1)^{k_1l}  h^*  [K_1,i(L)K_2] -\\
-(-1)^{k_1l}  i^h(h^*[K_1,L])K_2 + (-1)^{(k_1+l)k_2} 
    i^h(h^*[K_2,L])K_1.\endmultline$$
    
{\bf 3. } For $L\in \Om {l+1}^h$ and $K_i\in\Omh {k_i}$ we have
$$\multline i^h(L)(h^*  [K_1,K_2]) = h^*  [i^h(L)K_1,K_2]
    +(-1)^{k_1l}  h^*  [K_1,i^h(L)K_2]-\\
-(-1)^{k_1l}i^h(h^*[K_1,L])K_2 +(-1)^{(k_1+l)k_2} 
    i^h(h^*[K_2,L])K_1.\endmultline$$
    
{\bf 4. } For $K\in \Omh k$ and $L_i\in\Om {l_i+1}^h$ we have
$$\multline h\o [K,[L_1,L_2]^{\wedge,h}]\o \Lambda h =\\
= [h\o [K,L_1]\o\Lambda h,L_2]^{\wedge,h} + (-1)^{kl_1}  [L_1,h\o
    [K,L_2]\o \Lambda h]^{\wedge,h} -\\
-(-1)^{kl_1}h\o[i^h(L_1)K,L_2]\o\Lambda h -
    (-1)^{(l_1+k)l_2}h\o [i^h(L_2)K,L_1]\o\Lambda h.\endmultline$$
    
{\bf 5. } Finally for $K_i\in \Om {k_i}$ and $L_i\in
\Om {k_i+1}^h$ we have
$$\multline [ \Theta^h(K_1) + i^h(L_1),\Theta^h(K_2) + i^h(L_2)] =\\
= \Theta^h\bigl ([K_1,K_2]  + i(L_1)K_2 -(-1)^{k_1k_2}
      i(L_2)K_1\bigr )+\\
+i^h\bigl ( [L_1,L_2]^{\wedge,h} + (-1)^{k_2}  h\o
    [L_1,K_2]\o\Lambda h - (-1)^{k_1(k_2+1)}  h\o
    [L_,K_1]\o\Lambda h -\\
-h\o [K_1,\phi\o [K_2,\phi ]]\o\Lambda h +
    (-1)^{k_1k_2}  h\o [K_2, \phi\o [K_1,\phi ]]\o\Lambda h\bigr). 
\endmultline$$
\endproclaim

\demo{Proof: 1} 
$$\alignat2 [i^h(L),\Theta^h(K)] 
&= i^h(L) h^* \Theta (K)h^*  -(-1)^{kl}  h^* \Theta (K)  h^*  i^h(L) =\\
&= h^*\bigl (i(L)\Theta (K) -(-1)^{kl}\Theta(K)i(L)\bigr )h^*,
    &&\quad\text{ by 2.3.4.} \\
&=h^* \bigl (\Theta(i(L)K)+(-1)^ki([L,K])\bigr )  h^*,
    &&\quad\text{ by \cite{6,1.6}}\\
&=\Theta^h(i(L)K) + (-1)^k  i^h(h\o [L,K]\o\Lambda h),
    &&\text{ by 2.3.4 and 2.4.}\endalignat$$
    
{\bf 2. } $i^h(L)(h^*[K_1,K_2]) = h^*  i(L)[K_1,K_2]$ by 2.3.4
and the rest follows from \cite{6, 1.10}.

{\bf 3. } Note that for horizontal K we have $i(L)K=i^h(L)K$ and
use this in formula 2.

{\bf 4. } This follows from 1 by writing out the graded Jacobi
identity for the graded commutators, uses horizontality of $K$
and 2.3 several times. 

{\bf 5. } Collect all terms from 2.6.8 and 1.  \qed\enddemo

\subheading{2.9} The space of derivations over $h^*$ of  $ \Omega(M)$
is a graded module over the graded algebra $\Omega(M)$ with the
action $(\om\wedge D)\psi = \om\wedge D\psi$. The subspace $\text{
Der}^h\Omega(M)$ is stable under $\om\wedge\cdot$ if and only if $\om\in\Omeh
{}$.

\proclaim {Theorem} 1. For derivations $D_1,D_2$ over $h^*$ of
degree $k_1,k_2$, respectively, and
$\om\in\Omeh q$ we have $[\om\wedge D_1,D_2]=\om\wedge
[D_1,D_2]- (-1) ^{(q+k_1)k_2}D_2\om\wedge h^*D_1$.

2. For $L\in\Om {}$ and $\om\in\Ome {}$ we have $\om\wedge
i^h(L)=i^h(\om\wedge L)$.

3. For $K\in\Ome k$ and $\om\in\Ome q$ we have
$$(\om\wedge\Theta(K))\o h^* = \Theta^h(\om\wedge K)\o h^*
    +(-1)^{q+k-1}i^h(d\om\wedge (h\o K)).$$

4. For $K\in\Omh k$ and $\om\in\Omeh q$ we have
$$\om\wedge\Theta^h(K) = \Theta^h(\om\wedge K) +(-1)^{q+k-1}
    i^h(d^h\om\wedge (h\o K)).$$

5. $\Om {}^h$ is stable under multiplication with $\om\wedge.$
if and only if $\om\in\Omeh {}$. For $L_j\in\Om {l_j+1}^h$ and
$\om\in\Omeh q$ we have 
$$[\om\wedge L_1,L_2]^{\wedge,h} =
\om\wedge [L_1,L_2]^{\wedge,h}-
(-1)^{(q+l_1)l_2}\,i^h(L_2)\om\wedge(h\o L_1).$$
\endproclaim
The proof consists of straightforward computations.

\heading 3. Derivations on Principal Fiber Bundles and
Associated Bundles.\\ Liftings.\endheading

\subheading{3.1} Let $(E,p,M,S)$ be a fiber bundle and let $\phi
\in\Om 1$ be a connection for it as described in 1.1. The
cocurvature $\bar R$ is then zero. We consider the horizontal
lifting $\chi: TM\times_M E\to  TE$ and use it to define the mapping
$\chi_*:\Om k \to  \Omega^k(E;TE)$ by 
$$(\chi_*K)_u(\row X1k) := \chi\bigl (K(\row {T_up.X}1k),u\bigr).$$
Then $\chi_*K$ is horizontal with horizontal values: $h\o\chi_*K
= \chi_*K = \chi_*K\o\Lambda h$. For $\om\in\Ome q$ we have
$\chi_*(\om\wedge K) = p^*\om\wedge\chi_*K$, so $\chi_* : \Omega
(M;TM) \to  \Omega (E;TE)^h$ is a module homomorphism of degree
$0$ over the algebra homomorphism $p^*:\Omega (M)\to \Omega_{\text{hor}}(E)$.

\proclaim {Theorem} In this setting we have for 
$K,K_i\in\Omega(M;TM)$: 
\roster
\item $p^*\o i(K) = i(\chi_*K)\o p^* = i^h(\chi_*K)\o p^* 
    :\Omega (M) \to  \Omega_{\text{hor}}(E)$.
\item $p^*\o\Theta (K)  = \Theta(\chi_*K)\o p^* =
    \Theta^h(\chi_*K)\o p^* : \Omega (M) \to  \Omega\h (E)$.
\item $i(\chi_*K_1)\chi_*K_2 = i^h(\chi_*K_1)\chi_*K_2 =
    \chi_*(i(K_1)K_2)$.
\item The following is a homomorphism of graded Lie algebras:
    $$\chi_*:(\Om {},[\quad ,\quad]^\wedge)
    \longrightarrow (\Omega (E;TE)^h,[\quad,\quad]^{\wedge,h})
    \subset(\Omega (E;TE),[\quad,\quad]^\wedge)$$
\item  $\chi_*[K_1,K_2] = h\o [\chi_*K_1,\chi_*K_2] = h\o
    [\chi_*K_1,\chi_*K_2]\o \Lambda h$.
\endroster
\endproclaim
\demo{ Proof} \therosteritem1. $p^*\o i(K)$ and $i(\chi_*K)\o p^*$ 
are graded module derivations $:\Omega (M)\to  \Omega (E)$ over
the algebra homomorphism  
$p^*$ in a sense similar to 2.1:
 $p^*i(K)(\om\wedge\psi) =
p^*i(K)\om\wedge p^*\psi + (-1)^{(k-1)\vert
\om\vert} p^*\om\wedge p^*i(K)\psi$ and similarly for the other
expression. Both are zero on $C^\infty (M)$ and it suffices to
show that they agree on $\om\in\Ome 1$. This is an easy
computation. Furthermore $i^h(\chi_*K)\vert\Omega\h (E) =
i(\chi_*K)\vert \Omega\h (E)$, so the rest of 1 follows. 
\therosteritem2  follows from \therosteritem1 by expanding the
graded commutators. 
\therosteritem3.~Plug in the definitions and use 2.3.1.
\therosteritem4~follows from 2.3.3 and \cite{6,1.2.2}.
\therosteritem5~follows from 2 and some further considerations. 
\qed\enddemo 

\subheading{3.2} Let $(P,p,M,G)$ be a principal fiber bundle with
structure group $G$, and write $r:P\times G\to  P$ for the
principal right action. A connection $\phi^P:TP\to  VP$ in the
sense of 1.1 is a {\it principal connection} if and only if it is
$G$-equivariant i.e., $T(r^g)\o \phi^P = \phi^P\o T(r^g)$ for
all $g\in G$. Then $\phi^P_u = T_e(r_u).\varphi_u$, where
$\varphi\,$, a 1-form on $P$ with values in the Lie algebra of
$G$, is the usual description of a principal connection. Here we
used the convention $r(u,g)=r^g(u)=r_u(g)$ for $g\in G$ and $u\in
P$. The curvature $R$ of 1.3 corresponds to the negative of the
usual curvature $d\varphi + \tfrac12[\varphi,\varphi]$ and the
Bianchi identity of 1.4 corresponds to the usual Bianchi identity.

The $G$-equivariant graded derivations of $\Omega (P)$ are
exactly the $\Theta(K) + i(L)$  with $K$ and $L$ $G$-equivariant
i.e., $T(r^g)\o K = K\o\Lambda T(r^g)$ for all $g\in G$. This
follows from \cite{6, 2.1}. The $G$-equivariant derivations in
$\text{Der}^h\Omega (P)$ are exactly the $\Theta^h(K)+i^h(L)$
with $K\in \Omega\h (P;TP)$, $L\in\Omega  (P;TP)^h$ and $K$,$L$
$G$-equivariant. Theorem 3.1 can be applied and can be
complemented by $G$-equivariance.

\subheading{3.3} In the situation of 2.2 let us suppose furthermore,
that we have a smooth left action of the structure group $G$ on
a manifold $S$, $\ell :G\times S \to  S$. We consider the associated
bundle $(P[S]=P\times_GS,p,M,S)$ with its $G$-structure and
induced connection $\phi^{P[S]}$. Let $q:P\times S \to  P[S]$ be
the quotient mapping, which is also the projection of a principal
$G$-bundle. Then $\phi^{P[S]}  Tq  (X_u,Y_s) = Tq  (\phi^P 
X_u, Y_s)$ by \cite{7, 2.4}.

Now we want to analyze $q^*\o (\Theta^h(K)+i^h(L)): \Omega(P[S])
\to  \Omega(P\times S)$. For that we consider the associated
bundle $(P\times S = P\times_G(G\times S),p\o pr_1,M,G\times
S,G)$, where the structure group $G$ acts on $G\times S$ by left
translation on $G$ alone. Then the induced connection is
$\phi^{P\times S} = \phi^P\times Id_S: T(P\times S) = TP\times
TS \to  VP\times TS$  and we have the following 

\proclaim{Lemma} $Tq\o h^{P\times S} = h^{P[S]}\o Tq$ and
$D^{h^{P\times S}}\o q^* = q^*\o D^{h^{P[S]}}$ for the classical
covariant derivatives.
\endproclaim

The proof is obvious from the description of induced connections
given in \cite{7, 2.4}.

\proclaim{3.4. Theorem} In the situation of 3.3 we have:\newline
1. $q^*\o (h^{P[S]})^* = (h^{P\times S})^*\o q^*: \Omega
(P[S])  \to  \Omega  (P\times S)$.

2. For $K\in\Omega^k(P[S];T(P[S]))$ we have
$$\gather q^*\o i^{h^{P[S]}}(K) = i^{h^{P\times S}}(\chi^{P\times
S}{}_* K)\o q^*: \Omega (P[S])\to  \Omega\h (P\times S)\\
q^*\o \Theta^{h^{P[S]}}(K) = \Theta^{h^{P\times
S}}(\chi^{P\times S}{}_* K)\o q^*: \Omega (P[S]) \to  \Omega\h
(P\times S).\endgather$$
\endproclaim
\demo{ Proof} 1.  Use lemma 3.3 and the definition of $h^*$ in
2.2. For 2  use also 2.3.1 and $Tq\o \chi_*K = K\o\Lambda
Tq$. The last equation is easy.  \qed\enddemo 

\subheading{3.5} For a principal bundle $(P,p,M,G)$ let
$\rho:G\to GL(V)$ be a linear representation in a finite
dimensional vector space $V$. We consider the associated vector
bundle $(E=P\times_GV=P[V],p,M,V)$, a principal connection
$\phi^P$ on $P$ and the induced $G$-connection on $E$, which
gives rise to the covariant exterior derivative $\nabla^E \in
\text{Der}_1\Omega (M;E)$ as investigated in \cite{6, section
3}. We also have the mapping $q^\sharp :\Omega^k(M;E) \to 
\Omega^k(P,V)$, given by
$$(q^\sharp\Psi)_u(\row X1k) = q_u{}^{-1}(\Psi_{p(u)}(\row {T_up.X}1k).$$
Recall that $q_u:\{u\}\times V \to  E_{p(u)}$ is a linear
isomorphism. Then $q^\sharp$ is an isomorphism of $\Omega (M;E)$
onto the subspace $\Omega\h (P,V)^G$ of all horizontal and
$G$-equivariant forms.

\proclaim {Theorem} 1. $q^\sharp\o\nabla^E = D^{h^P}\o q^\sharp =
d^{h^P}\o q^\sharp$.

2. For $K\in \Om k$ we have:
$$\gather q^\sharp\o i(K) = i(\chi^P{}_*K)\o
    q^\sharp = i^{h^P}(\chi^P{}_*K)\o q^\sharp\\
q^\sharp\o\Theta_{\nabla^E}(K) = \Theta^{h^P}(\chi^P{}_*K)\o
    q^\sharp : \Omega (M;E) \to  \Omega\h (P,V)^G.\endgather$$  
\endproclaim
\demo{ Proof} 1. A proof of this is buried in \cite{3, p 76,
p 115}. A global proof is possible, but it needs a more detailed
description of the passage from $\phi^P$ to $\nabla^E$, e.g. the
connector $K:TE \to  E$. We will not go into that here.

2. $q^\sharp\o i(K)$ and $i(\chi^P{}_*P)\o
q^\sharp$ are both derivations from the $\omega (M)$-module
$\Omega (M;E)$ into the $\Omega (P)$-module $\Omega (P,V)$ over
the algebra homomorphism $q^*:\Omega (M) \to  \Omega (P)$. Both
vanish on 0-forms. Thus it suffices to check that they coincide
on $\Psi\in\Omega^1(M;E)$, which may be done by plugging in the
definitions. $i(\chi^P{}_*K)$ and $i^{h^P}(\chi^P{}_*K)$
coincide on horizontal forms, so the first assertion follows.
The second assertion follows from the first one, lemma 3.3, and
\cite{6, 3.7}. Note finally that $q^\sharp\o\Theta_{\nabla^E}(K)
\not= \Theta(\chi^P{}_*K)\o q^\sharp$ in general. \qed\enddemo 

\subheading{3.6} In the setting  of 3.6, for $\Xi \in
\Omega^k(M;L(E,E))$, let $q^\sharp \Xi \in \Omega (P;L(V,V))$ be
given by  
$$(q^\sharp\Xi)_u(\row X1k) =
q_u{}^{-1}\o\Xi_{p(u)}(\row {T_up.X}1k)\o q_u:V \to  V,$$  
where $L(E,E)=P[L(V,V),Ad\o\rho] $ and $q^\sharp\Xi$ defined here
coincides with that of 3.5.

Furthermore recall from \cite{6, section 3} that any graded $\Omega
(M)$-module homomorphism $\Omega (M;E) \to  \Omega (M;E)$ of
degree $k$ is of the 
form $\mu (\Xi)$ for unique $\Xi \in \Omega^k(M;L(E,E))$, where 
$$(\mu (\Xi)\Psi)(\row X1{k+q}) = \tfrac1{k!\, q!}
\sum_\sigma\varepsilon(\sigma).\Xi(\row X{\sigma 1}{\sigma k}).
\Psi(\row X{\sigma(k+1)}{\sigma(k+q)}).$$ 
Then any graded
derivation D of the $\Omega (M)$-module $\Omega (M;E)$ can be
uniquely written in the form $D= \Theta_{\nabla^E}(K)+i(L)+\mu(\Xi)$.

\proclaim{Theorem} For $\Xi\in\Omega  (M;L(E,E))$ we have 
$$q^\sharp\o\mu(\Xi) = \mu(q^\sharp\Xi)\o q^\sharp =
\mu^{h^P}(q^\sharp\Xi)\o q^\sharp:\Omega (M;E) \to  \Omega\h
(P,V)^G,$$ 
where for $\xi\in \Omega^k(P,L(V,V))$ the module
homomorphism $\mu(\xi)$ of $\Omega (P,V)$ is as above for the
trivial vector bundle $P\times V \to  P$ and $\mu^h(\xi)$ is
given by  
$$\multline (\mu^h(\xi)\om)(\row X1{k+q}) = \\
= \frac1{k!\,q!} \sum_{\sigma\in\Cal S_{k+q}}
    \varepsilon(\sigma)\xi(\row X{\sigma 1}
    {\sigma k})\om(\row {h^P\,X}{\sigma(k+1)}{\sigma(k+q)}).
\endmultline$$
\endproclaim

The proof is straightforward.

\enddocument